\documentclass[12pt]{article}

\usepackage{a4wide,amsfonts,amsmath,amssymb}

\def\i{{\bf i}}
\def\j{{\bf j}}

\def\k{{\bf k}}

\def\A{{\bf A}}

\def\F{{\bf F}}

\def\H{{\bf H}}

\def\A{\vec{\bf A}}

\def\x{{\bf x}}

\def\F {{\vec{\bf F}}}

\def\i{{\bf i}}

\def\k{{\bf k}}

\def\rot{{\operatorname{rot}}}

\def\B{\vec{\bf B}}

\begin{document}

	\title{Invariants of magnetic lines for Yang-Mills solutions}
	\author{Akhmet'ev P.M.\footnote{National Research University “Higher School of Economics”, Moscow, Russia;
			IZMIRAN, Troitsk, Moscow, Russia} and Dvornikov M.S.\footnote{IZMIRAN, Troitsk, Moscow, Russia}}
	\maketitle

\begin{abstract}
We construct a new Yang-Mills $3D$-solution on the space of negative scalar curvarure. We discuss a problem of 	non-abelian gauge symmetry is broken with the assumption that a scalar curvature of the domain is a negative small parameter.
	In this case we use the following fact: a geometrical scale related with Vassiliev's discriminant of magnetic lines coincids with a physical Kolmogorov scale. This gives
an estimation of $\alpha$-effect by the dispersion of the asymptotic ergodic Hopf invariant in the limite with a negative scalar curvature parameter.
\end{abstract}

\section{Introduction}

Since the theoretical discovery of the Yang--Mills theory in 1954 in Ref.~\cite{YanMil54}, the importance of non-abelian gauge fields for modern physics is enormous. It is mainly because of the fact that the experimentally confirmed standard model of elementary particles is based on $U(1)\times SU(2)\times SU(3)$ gauge group. The main problem in the study of non-abelian gauge fields is that the corresponding wave equations are generically nonlinear. It leads to the appearance of various non-perturbative solutions, like instantons, which play an important role in the construction of the models for the quantum chromodynamics (QCD) vacuum~\cite{ShuSch97}. Numerous non-perturbative methods for the study of Yang--Mills fields are reviewed in Ref.~\cite{Hub20}.

Despite the fact that the non-abelian gauge symmetry of the standard model is broken presently owing to the Higgs mechanism, this symmetry was restored in the early universe before the electroweak phase transition (EWPT) which happened at $T = 100\,\text{GeV}$ or $t \sim 10^{-10}\,\text{s}$ since the Big Bang~\cite{GorRub11}. All particles, including gauge bosons, were massless before EWPT. The behavior of gauge bosons before EWPT can influence the evolution of Maxwell electromagnetic fields in the present universe. It is applied especially to the hypercharge field which is the progenitor of photons~\cite{DvoSem13}. If a nonzero hypermagnetic field existed in the early universe, it can result in the leptogenesis and, hence, the baryogenesis owing to the analogue of the Adler anomaly~\cite{GioSha98}.

Non-abelian gauge fields in the early universe evolve in the non-Minkowski spacetime. The spacetime of the expanding universe has the Friedmann--Robertson--Walker (FRW) metric. Despite the fact that FRW metric is conformally flat, the influence of nontrivial geometry on the evolution of a Yang--Mills field has to be studied. Various quantum effects involving Yang--Mills fields in curved spacetime are described in Ref.~\cite{GriMos80}. The quantization of Yang--Mills fields in curved spacetime was reviewed in Ref.~\cite{Hol08}.

The powerful tool for the study of the electromagnetic fields evolution in highly conducting medium involves the magnetic hydrodynamics (MHD) approximation. In this approach, one can exclude the electric field and deal with a magnetic field only. As a result one obtains the induction equation~\cite{LanLif84}. If P-symmetry is violated in a system, the induction equation acquires the $\alpha$-term responsible for the magnetic field instability. The mechanism for the $\alpha$-term generation in solar physics is described in Ref.~\cite{Sti04}.

However, if one deals with the electromagnetic field in a curved spacetime, the induction equation is difficult to formulate since one cannot separate the electromagnetic field tensor into electric and magnetic fields in an invariant manner. One can make it in the FRW metric using comoving variables (see, e.g., Ref.~\cite{DvoSem14}). However, in this situation, one deals with the spectra of the magnetic energy and the magnetic helicity. The same problem arises for non-abelian fields. In this case, the magnetic field is a matrix valued function. In the present work, we consider one of the ways to introduce the $\alpha$-parameter for Yang--Mills fields in the invariant manner in a curved spacetime.

It is clear that properties of $SU(2)$ cannot be totally described by means of helicity integral
and higher invariant of magnetic lines are into account. A starting point for the present paper is
a work \cite{H-M}, where the authors shows that Chern-Simons integral is related with higher Massey product of magnetic lines of Yang-Milles solution. However, the Massey product are not universal invariants in MHD this invariant do not satisfy ergodic asymptotic properties, as the helicity integral. In the textbook \cite{A-Kh}, based on an Arnold's paper \cite{Arn}, the asymptotic ergodic Hopf invariant is introduced, this invariant is an expression of the helicity density of magnetic lines in an invariant form. Several authors \cite{R-Y-H-W} introduced a notion of field line helicity, which coincides with the convolution of asymptotic Hopf invariant in the case when magnetic line are concentrated in a closed domain and are tangent to the boundary of the domain. But in a general case, when magnetic lines are across to a component of the boundary of the domain, an extra terms in the construction of the field line helicity is required. At present a various modifications of helicity integrals are studied: relative magnetic helicity \cite{T-V}, winding magnetic helicity \cite{X-P-Y}. All this constructions are applied to Maxwellian magnetic field, but allow no satisfactory description of properties of $SU(2)$-fields (see Section~2 for helicity term in CS-integral).

In the paper by the first author \cite{A1} a construction of higher invariant of magnetic lines for Maxwell field are introduced in an asymptotic and ergodic form and in \cite{A2} are related with  the $\beta$-term  in the mean field theory. From mathematical point of view there are difficulties in the definition. Even in the case of helicity, the conception of the asymptotic ergodic Hopf invariant in the case magnetic fields is not generic  gives difficulties, see \cite{A-Kh} Ch.III, Sec. 4.3. In the case magnetic lines are closed (this case is also non-generic, one speaks about magnetic field concentrated in  magnetic tubes) we get well-defined notion of density of  higher magnetic invariant ($M$-invariant), but not of higher Massey invariant originally was introduced in MHD in \cite{M-R}. 

As we have shown linking number of field lines, or, the magnetic helicity, allows no description of
magnetic fields in a curved domain, higher invariants of magnetic lines are required. Based on Kolmogorov's approach on hydrodynamics, which is known in MHD as the Kolmogorov's magnetic spectrum, we study a contribution of the scalar curvature parameter of the domain to $\alpha$-term in the mean field equation \cite{M}. The construction is based on the idea that a scale of the space, which a priori give contribution to higher invariant of magnetic lines, has to lead invariant gauge. This is observed in the case the density of magnetic lines is distributed as a magnetic flow of Kolmogorov spectrum \cite{A2}. If we assume that $M_5$-invariant of magnetic lines distributes in the domain uniformly with the scale, we get the Kolmogorov's $\beta$-term scale distribution. But, the distribution of magnetic lines with uniformly $M_5$-density depends on the scalar curvature parameter and its asymptotic for small curvature parameter is easily calculated. With an additional assumption that $\alpha$-term is distributed with the same scale, we get the required estimation of the asymptotic of the $\alpha$-parameter with respect to the scalar curvature of the domain.

This work is organized in the following way. In Sec.~1, we 
construct a Yang-Milles solution on the standard 3D-sphere $S^3$ with positive scalar curvature and give a short history of the question. In Sec.~2, we generalize  the solution for  the domain with negative scalar curvature. We use calculations from the textbook \cite{A-Kh} and
give a Yang-Milles solution for this domain.
The Yang-Milles solution is bracken, and we get a Maxwell field with the Kolmogorov distribution, which is closed  to a component of Yang-Milles field. We have a non-compact domain in this case, and we assume that
a field is quasi-periodic. This assumption means, that the Chern-Simons wave, which is described  the large scale component of magnetic field,
is closed to geodesic flow in the domain. In Sec.~3 we recall properties of $M_5^{sim}$ invariant \cite{A2} and calculate this invariant for a component of Yang-Milles field first, and for closed Maxwell magnetic field, in geometrical terms of Vassiliev discriminant theory \cite{Arn2}.
This geometrical approach gives an integral contribution of asymptotic of linking number of closed magnetic lines. Physically this interpreted as linking numbers of magnetic lines in the magnetic flow in mean-field equations, and this gives an asymptotic of the $\alpha$-term by the scalar curvature parameter.  In Sec.~4 we give conclusions.

\section{Yang-Mills field and Chern-Simons invariant}

In this section a connection $\A^R$ in $SU(2)$-bundle $\eta^R$ over the standard 3D sphere $S^3$ is studied. The bundle $\eta^R$ is isomorphic to the trivial bundle $S^3 \times \H$. 
Let us consider $S^3$ as the unite sphere in $\H$ with the base vectors  $\{+1,+\i,+\j,+\k\}$. We identify the bundle $\eta^R$ by the isomorphism: 
$$ \phi^R: \H \to \eta^R, \quad \phi^R(\xi) = \xi \circ \x, \quad \x \in S^3, \quad \xi \in \H. $$
In this formula $\circ$ is the quaternion product. This formula determines $\A^R$ as $\phi^R$-image of imaginary quaternions base vectors $\{\i,\j,\k\}$. By this definition $\A^R$ is the Levi-Civita connection over the bundle $\eta^R$ (as the real bundle this is 4D bundle). In the construction we
use the metric, and we write-down connections and its curvatures as vectors in the standard coordinate system.  There is an isomorphism  $T(S^3)= Im(\eta^R)$, where $T(S^3)$ is the tangent bundle, $Im(\eta^R) \subset \eta^R$ is a 3D subbundle generated by imaginary quaternions. Also this formula determines the right-polarized frame on $T(S^3)$ by the $\phi^R$-image of the imaginary base quaternions.  

The right-action by the base quaternions in the fiber of $\eta^R$ is an isometry of the fiber, this isometry is described by the corresponding matrix, which is called Pauli matrix.

Yang-Milles field is a $2$-form $\F$, which is the curvature of the connection $\A$ in $SU(2)$-bundle $\eta$. This form satisfies the equations:
\begin{eqnarray}\label{YM}
	D\F=0; \quad D^*(\F) = 0.
\end{eqnarray}
The curvature $\F^R$ of the connection $\A^R$ satisfies the equation (\ref{YM}(1)), but not the equation (\ref{YM}(2)). To get a Yang-Milles field, we replace $\A^R \mapsto \A^L$ as following.

Recall, the fiber $\eta_{\x}$ over a point $\x \in S^3$ is defined by the right action on the quaternion $\x$ on the standard fiber $\H$. Let us define the left action  by the conjugated quaternion $\bar{\x}$, namely:
$$ \phi^L:  \H \to \eta^L, \quad \phi^L(\xi) = \bar{\x} \circ \xi, \quad \x \in S^3, \quad \xi \in \H; $$
$$ \bar{\x}= \overline{a+b\i+c\j+d\k} = a -b\i - c\j -d\k. $$

Let us define the connection $\A^L$ in $\eta^L$.
Let us define the Riemannian metric on $\eta^L$ analogously with the case $\eta^R$ above.
The fiber 
$\eta^L_{\x}$, $\x=+1 \in S^3,$ is isometric to $\H$. The triple
$$ \A^L = (\A_i^L,\A_j^L,\A_k^L) $$
is defined as the result of the left translation of the base quaternions $\i_{\eta},\j_{\eta},\k_{\eta}\}$ from the fiber $\eta^L_{+1}$ into the space of the bundle.
Because $Im(\eta^L)$ is isomorphic to $T(S^3)$, we may assume that $\A^L$ is the frame in $T(S^3)$,
which is a $1$-form of the connection.  
In the prescribed coordinate system we may consider $\A^L$ as vectors, by definition 
$(\A_i^L,\A_j^L,\A_k^L) = (\bar{\i}, \bar{\j}, \bar{\k})$,
where in the right hand side we get the conjugated quaternion frame on $T(S^3)$. 

Let us define a $2$-form $\F^L$, the Yang-Milles field, using the Riemannian metric, this form is defined as the vector-field by the formula:
$$ \F^L = (-\A_i^L,-\A_j^L,-\A_k^L). $$
In this case the equation  (\ref{YM};1) is written as: 
$$ \F^L = \rot \ \A_l^L + [\A^L_{l+1},\A^L_{l+2}], $$
where $[,]$ is the commutator  $[\A_l^L,\A_{l+1}^L]=\A_{l+2}^L$, $l=\i,\j,\k$, $\rot$
is the standard vorticity operator on $S^3$. 
Because
$$D^{\ast}d(\A^L)=\rot \ \rot(\A_l^L)+[\rot(\A_{l+1}^L),\A_{l+2}^L]+[\A_{l+1}^L,\rot(\A_{l+2}^L)]=0,$$
the equation (\ref{YM}(2)) is written in the form: 
$$ \rot[\A_{l+1}^L,\A_{l+2}^L] + [\A_{l+1}^L,[\A_{l}^L,\A_{l+1}^L]] + [\A_{l+2}^L,[\A_{l}^L,\A_{l+2}^L]] =0. $$
The equations  (\ref{YM}(1))(\ref{YM}(2)) are satisfied because of the equation: 
$$ \rot \A_l^L = -2 \A_l^L, \quad [\A_l^L,\A_{l+1}^L]=\A_{l+2}^L. $$

The following integral:
\begin{eqnarray}\label{CS}
	CS(\F^L)= \int_{S^3}  tr(\A^L \wedge d\A^L + \A^L \wedge \A^L \wedge \A^L) 
\end{eqnarray}
is called the Chern-Simons functional.

The first term of the integral (\ref{CS}) is defined as the average linking number of magnetic lines of the fields:
\begin{eqnarray}\label{AB}
	\B_{\i}^L = \rot\A_{\i}^L,  \B_{\j}^L = \rot \ \A_{\j}^L, \B_{\k}^L = \rot\A_{\k}^L .
\end{eqnarray}
The second term is an average number of higher linking numbers of magnetic lines of the triple:  $\{\B_{\i}^L, \B_{\j}^L, \B_{\k}^L\}$.

\subsection*{An history of the question}
In a paper by Hornig and Mayer \cite{H-M} for the right-hand fields the equation is satisfied:
$$ \rot \A_l^R = +2 \A_l^R. $$
In this case Yang-Milles equation (\ref{YM}(2)) is not satisfied.

By the quaternion conjugation $\A^R \mapsto \A^{L}$
the proper number of the operator $\rot$ changes the sign.
The equation:
$$[\A_l^{L},\A_{l+1}^{L}]=\A_{l+2}^{L} $$
is satisfied, because the base 
$\eta^L_{+1}$ remains standard.

The space $S^3$ from the holomorphic point of view is the right-polarized space. One  may considered this space alternatively as the left-polarized space by a different representation. 
In this construction $S^3$ is  the $Spin$-covering for the geodesic flow on the standard $2D$-sphere  $S^2$.
Let us shows that in this caser Yang-Mills equations are generalized.

\section{Geodesic flows of positive and negative curvature }

In this section we use definition and calculations from  \cite{A-Kh} Ch.5; 4.4.
The previous calculation for Yang-Mills field on $S^3$  are reformulated for the $Spin$-covering
over the space
$$ SO(3) = \{\xi \in T_{\ast}(S^2): \vert\vert \xi \vert \vert = 1 \}. $$
All the calculations remains analogous, if we replace $S^2$ into the Lobachevski plane
$\Lambda$.

Let us consider the space
$$ T\Lambda = \{\xi \in T_{\ast}(\Lambda): \vert\vert \xi \vert \vert = 1 \}, $$
where $\Lambda$ is the Lobachevski  plane with the unite negative scalar curvature $-1$.
Let us define the right triple of vector fields:
$$\{ f^R,e^R,\tilde{e}^R \}.$$
The field $f^R$ is the isometric rotation trough a positive angle in the prescribed fiber,
$e^R$ is the geodesic flow,  $\tilde{e}^R$ is the conjugated geodesic flow on the Lobachevski  plane. 
The sectional curvature on the plane
$(e^R,\tilde{e}^R)$ equals to $-1$, sectional curvatures in the planes
$(f^R,e^R)$, $(f^R,\tilde{e}^R)$ equal to $-2$. The last two equations follow is proved as following.
Replace $\Lambda$ by $S^2$, the scale in the vertical fiber $S^1 \subset S^3$ over $\tilde{T}(S^2)$ is normalized by the standard metric in the base of the double covering $\times 2: S^1 \to S^1$. The calculation for $S^2 \mapsto \Lambda$ is analogous, the only difference that the scalar positive curvature in all sections changes the sign.

Take the parameter $\lambda > 0$, the scale of the metric of the vertical coordinate
$f^R$ such that the following equations are satisfied:
$$ (e^R,e^R)=(\tilde{e}^R,\tilde{e}^R)=1; \quad (f^R,f^R)=\lambda^2. $$
In this case the following equations are satisfied:
$$ \rot \ e^R = -\frac{e^R}{\lambda}; \quad \rot \ \tilde{e^R} = - \frac{\tilde{e^R}}{\lambda}; \quad \rot \ f^R = f^R. $$ 
The volume form, generated by this vector fields equals to
$\mu(e^R,\tilde{e}^R,f^R)=\lambda$.
The parameter
$\lambda$, by this equation, changed the phase volume into $\lambda$ times. 

Let us pass into the left base, which is given by the following triple of the vector fields: 
$$\{ f^L=- \frac{f^R}{\lambda}, \ e^L=e^R, \ \tilde{e}^L=\tilde{e}^R \}.$$
For this base we get the relation: 
$$ \rot \ e^L = -\frac{e^L}{\lambda}; \quad \rot \ \tilde{e^L} = - \frac{\tilde{e^L}}{\lambda}; \quad \rot \ f^L = - \frac{f^L}{\lambda}. $$

Let us denote by $\tilde{T}(\Lambda)$ the double covering over $T(\Lambda)$. The Riemannian metric
in this space is taken analogously with the case
$\tilde{T}(S^2)$: in the horizontal plane this metric relates with  the metric on $T(\Lambda)$, 
in the vertical circle the length of the vector is twice greater then the angle length on the base. For  $\lambda=1$ we get a manifold with a constant scalar curvature, because all sectional curvatures equal $-1$. In a general case the sectional curvature in the horizontal plane $(e^L,\tilde{e}^L)$ equals to $-1$, and sectional curvatures in the planes  $(f^L,e^L)$, $(f^L,\tilde{e}^L)$ equal to  $ \lambda^{-1}$.

In the Riemannian manifold
$\tilde{T}(\Lambda)$ the base fields is re-denoted by the following gauge:  
$$f^L \mapsto -2f^L, \quad e^L \mapsto e^L, \quad \tilde{e}^L \mapsto \tilde{e}^L,$$ 
the resulting fields are not re-denoted.

For the fields on
$\tilde{T}(\Lambda)$ the following relations are satisfied: 
$$ \rot \ e^L = -2\frac{e^L}{\lambda}; \quad \rot \ \tilde{e}^L = -2 \frac{\tilde{e}^L}{\lambda}; \quad \rot \ f^L = -2 \frac{ f^L}{\lambda}. $$ 
Let us consider the case
$\lambda=1$. As in the previous section, the triple $\{ f^L,e^L,\tilde{e}^L \}$
determines the Yang-Milles solution.

Let us consider the case of an arbitrary
$\lambda > 0$. The scale of the space is normalized by its volume.  By this scale the fields 
\begin{eqnarray}\label{A}
	\A_{\i}^{\lambda}=  \lambda^{\frac{1}{3}}f^L;\quad \A_{\j}^{\lambda}= \lambda^{\frac{1}{3}} e^L, \quad \A_{\k}^{\lambda}=\lambda^{\frac{1}{3}} \tilde{e}^L 
\end{eqnarray}
are normalized base fields. The sectional curvature in the plane $(\j,\k)$ (the horizontal plane)  equals to $\lambda^{-\frac{2}{3}}$. The sectional curvature in the plane:
$(\i,\j)$, $(\i,\k)$ contains $\lambda^{-\frac{5}{3}}$. The following equations are satisfied:
$$ \rot \A^{\lambda}_{\i} = - \frac{2}{\lambda} \A^{\lambda}_{\i}, \quad \rot \A^{\lambda}_{\j} = -\frac{2}{\lambda} \A^{\lambda}_{\j}, \quad \rot \A^{\lambda}_{\k} = -\frac{2}{\lambda} \A^{\lambda}_{\k}. $$ 
In the case  $\lambda \to +\infty$ the space become flat.  

Because the scale of the space determines a proportional scale of the brackets,
the relation
$$[\A_{\i}^{\lambda},\A_{\j}^{\lambda}]=\A_k^{\lambda}, \quad [\A_{\k}^{\lambda},\A_{\i}^{\lambda}]=\A_{\j}^{\lambda}, \quad [\A_{\j}^{\lambda},\A_{\k}^{\lambda}]=\A_{\i}^{\lambda}$$ 
remains as in the case
$\lambda=1$.

Let us denote $\B_{\i}^{\lambda} =  -\frac{2}{\lambda} \A_{\i}^{\lambda}$, $\B_{\j}^{\lambda} =  -\frac{2}{\lambda} \A_{\j}^{\lambda}$, $\B_{\k}^{\lambda} =  -\frac{2}{\lambda} \A_{\k}^{\lambda}$. Let us calculate densities of magnetic helicities. We get:
$$ (\A_{\i}^{\lambda},\B_{\i}^{\lambda}) = 
(\A_{\j}^{\lambda},\B_{\j}^{\lambda}) = 
(\A_{\k}^{\lambda},\B_{\k}^{\lambda}) = -2. $$ 
We see that a sub-integral kernel of the each field is constant and depends not of the parameter
$\lambda$. As the result, the density of magnetic helicities, the first term in the integral,
depends not on sectional curvatures.

Let us consider the second term of the integral (\ref{CS}) as a function of the parameter $\lambda$.
We get that the sub-integral 3-form
\begin{eqnarray}\label{form}
	\A_{\i}^{\lambda} \wedge \A_{\j}^{\lambda} \wedge \A_{\k}^{\lambda}  
\end{eqnarray} 
is proportional to  $\lambda^{-1}$. Therefore the density of the second term 
in \ref{CS} is proportional to  $\lambda^{-1}$.

Let us present an easy geometrical reason of this calculation. 
Accordingly Hornig and Mayer \cite{H-M},  an analogous case of magnetic fields on
$\tilde{T}(S^2)$), the density $\rho$ of the second term in the integral  (\ref{CS})
is given by the triples magnetic lines distribution, which are detected by its starting points and its tangent vectors, and which are pairwise linked in the phase space. 

Let us consider a point
$x \in \tilde{T}\Lambda$, consider its projection  $\bar{x}$ on $\Lambda$ and 
calculate  the density of magnetic lines, starting in a neighborhood of $\bar{x}$, which are pairwise intersected,
as a function of the parameter $\lambda$.

Let us consider $3$ concentric circles of a great radius
$R_1,R_2,R_3$ in the central point $\bar{x}$. The ratio of the parallelism angle from 
$x_1 \in R_1$ to the line $l_2$  and the perimeter of the circle  $R_1$
equals to $\lambda^{-\frac{1}{3}}$, when 
$R_1 \to + \infty$. The integral asymptotic average of the geodesic triangles 
when radii tends to $+\infty$ is of the order 
$\lambda^{-1}$, because pairwise intersections are independently. As the result, we get
$\rho \sim \lambda^{-1}$, this gives the asymptotic of the second term in (\ref{CS}), as a function of the scalar curvature of the horizontal section.

\section{$M$-invariant of magnetic lines\label{sec:CURV}}

We assume that Yang-Mills equation for 
(\ref{A}), is broken.  Let us consider the only potential, say,  the potential $\A_j$, which generates the magnetic field  $\B_{\j}$, which is closed to the Yang-Milles component.
Let us estimate the square of the density distribution of triangles of magnetic lines  of the field $\B_{\j}$ on the horizontal plane. This  value is of the order $\sim \lambda^{-\frac{2}{3}}$.
An analogous calculation is presented in 
\cite{A3}, Lemma 5.1, with the assumption that there exist an additional symmetry of magnetic lines. (Note, that in \cite{A3} a scale of the angle and, therefore, the Riemannian metric in $\tilde{T}(\Lambda)$ is taken an alternative way, this leads to non-suitable choice of the parameter $a$ in the statement of Main Theorem in Section $3$). 

To estimate the square of the  density, let us take all  triangles, with intersections of the sides  
on angles $\ge \varepsilon$, where
$\varepsilon > 0$ is a small positive parameter. Then, let us take the limit  $\varepsilon \to 0+$. This limit is of the order $\sim \lambda^{-\frac{1}{3}}$, the square is of the order 
$\sim \lambda^{-\frac{2}{3}}$

Let us consider the rescaling into $l$-times: 
$\x \mapsto l\x$. In this case the volume form is changed in $l^3$ times. The density of the second term in the integral (\ref{CS}) remains unchanged, because each potential (\ref{A}), which is included in  (\ref{form}) is changed by the factor $l^{-1}$.
In the mean field theory \cite{M}, the magnetic field
$\B_{\j}$ can be consider as a large scale magnetic field. From physical assumption that a small-scale magnetic field is Kolmogorov distributed, with the exponent   $-\frac{5}{3}$.
By this assumption, the turbulent diffusion generate a magnetic flow with the exponent  $-\frac{7}{6}$ (the $\beta$-term in the mean field equation). Moreover, because the $\alpha$-term, the small-scale component twists a little magnetic lines of the large-scale magnetic field  $\B_{\j}$. A origin of the $\alpha$-effect and its amplitude is outside of our investigation. We only assume that $\alpha$ and $\beta$-terms in the mean field equation corresponds each other, and the $\alpha$-term is small. With this assumption we estimate the $\alpha$-term with respect to the curvature parameter $\lambda$. 

Accordingly to 
\cite{A1}, a higher asymptotic ergodic invariant  $M^{sim}_5$ of magnetic lines has a density exponent,
which coincides with the exponent of the magnetic flow of the Kolmogorov spectrum in $k$-space.  
In the  density of
$M^{sim}_5$-invariant for the magnetic field $\B_{\j}$ is presented an only non-trivial term. This is the main term described in Theorem  17 \cite{A2} (an auxiliary formula), with respect to the contribution of linking numbers of magnetic lines, which is presented in Theorem
17 \cite{A1}. This term is quadratic with respect to the magnetic flow, it calculates the square of the density of geodesic triangles of the projections of Chern-Simons magnetic line. 
Each triangle is taken proportional to the product of squares of magnetic flows.
Accordingly to the combinatorial formula
of $M_5^{sim}$-invariant\cite{A1}, in a random quintuple of magnetic lines all triples of lines are considered. In the case a triple represents a triangle, its contribution is taken proportionally to the product of $10$ pairwise linking coefficients in the considered quintuple of magnetic lines.

We assume that the parameter
$\lambda^{-1}$ is small. Let us divide quintuples of magnetic lines into two groups:
quintuples of pairwise far lines and all other quintuples, which contain at least two closed magnetic lines.  A pair of far magnetic lines, which represents by a pair of curves in the horizontal plane, closed to geodesics, are linked with the coefficient $-1$ (the only intersection points of this pair), or, unliked, when the projection are disjoint. Below it was shown that pairwise coefficients are distributed with a density which is not depended on $\lambda$,
the average values of the square of triangles is of the order
$\lambda^{-\frac{2}{3}}$. Because we assume that the integral over closed magnetic lines is small, we may assume that a contribution of triangles with at least two closed sides into $M^{sim}_5$-integral is proportional to  
$\alpha^2$, because this term determines a contribution of the second moment of the density of Hopf asymptotic ergodic invariant. The dispersion of the Hopf asymptotic ergodic invariant is investigated in \cite{A-V}. From this we get that an assumption that  $\alpha$ and $\beta$-terms in the mean field equation are balanced, we get the asymptotic $\alpha \sim \lambda^{-\frac{1}{3}}$.

\section{Discussion}
{\bf The well-known Gauss invariant (and its modern generalization suggested
	by Arnold) plays an important role in magnetic field generation in 
	celestial bodies. The Arnold form of the Gauss integral admits a gauge invariant large-scale density. One may assume that the second-order  time-derivative of this density is homogenous in the domain by a small-scale velocity, this gives an additional property for large-scale magnetic flow.}
In the present work, using this idea we have studied the following question: if the $\alpha$-parameter in dynamo mean field equation depends on the scalar curvature of the space.

We assume that the large-scale magnetic field, Chern-Simons magnetic wave, is closed to a component of a Yang-Mills solution and a small scale magnetic field satisfies the Kolmogorov distribution in $k$-space. Because of this assumption, the $\beta$-term is predicted by singularities of field lines of large scale magnetic field, the $\alpha$-term and $\beta$-term have to get a common distribution. 
This gives an asymptotic of $\alpha$-term with respect to the scalar curvature parameter.

After a short introduction in Sec.~1, we have constructed a new solution of Yang-Mills equation on the standard 3D sphere, based on an ansatz proposed by Hornig and Mayer \cite{H-M}, in Sec.~2. In Sec.~3, we have translated Yang-Mills solutions for the space of a negative scalar curvature, Chern-Simons wave is constructed as geodesic flow on the Lobachevsky plane. We have used definition and calculations from Arnold-Khesin textbook \cite{A-Kh}. 

In Sec.~4 we have proved the main result, assuming that the Chern-Simons wave, which is closed to a component of Yang-Milles solution, represents a large scale magnetic Maxwell field and admits the Kolmogorov's distribution. This assumption with the main result of the paper \cite{A1} by the first author proves that the magnetic flow is visualized using magnetic lines singularities (using geodesic triangles) in the framework of Vassiliev theory. 
Because calculations of singularities of fields lines depend on the curvature in the domain, we get a relationship of physical parameters in mean fields equation with the curvature. This gives an asymptotic for a small curvature parameter. 

{\bf Our results in Sec.~\ref{sec:CURV} require that the curvature of a spacetime to be nonzero, namely, negative. According to the measurements present in Ref.~\cite{Agh20}, our universe is flat with the high level of accuracy. However, one can consider perturbations of a metric which can result in a nonzero curvature.}

The contribution of the gravity to the dynamo $\alpha$-parameter can result from a nonzero torsion terms (see, e.g., Ref.~\cite{Bam12}). In the present work, we have studied another mechanism for the contribution of gravity to the $\alpha$-parameter in torsionless spaces. It is valid for non-abelian Yang-Mills fields. This contribution can be important, e.g., for the study of gravity perturbations on the evolution of a hypercharge field in the early universe.

The perturbations of gravitational field are known to cause the deviation of the beam of light propagation from the geodesics. {\bf For the first time, this phenomenon was predicted in Ref.~\cite{Zel64}. This subject was reviewed in Ref.~\cite{Dol17}.} The recent mathematical development of this problem is given in Ref.~\cite{Sok23}. Thus, observing the propagation of photons, emitted by distant sources, one can study the spectral distribution of density perturbations in the universe~\cite{Fle19}.

According to the Einstein equations, matter density perturbations result in fluctuations of gravitational fields. Relying on our results, one can study  the contribution of these density fluctuations to the $\alpha$-parameter. In this vein, we can explore the dynamics of stochastic primordial (hyper-)magnetic fields. Note that the research of stochastic cosmic magnetic fields was also carried out in Ref.~\cite{Mac02}. The proposed study is important for observational constraints on primordial magnetic fields~\cite{Ade16}.

\section*{Acknowledgments}

We are thankful to V.B.Semikoz for useful discussions.
\[  \]

\end{document}